\documentclass[a4paper,12pt]{article}
\begin{document}

\author{S.V. Ludkowski.}
\title{Topological spaces related with topological quasigroups.}
\date{18 December 2023.}
\maketitle

\begin{abstract}
The article is devoted to a structure of topological spaces related
with topological quasigroups. Regular and complete spaces over
topological quasigroups are studied. Separations and embeddings are
also investigated for them. Their homeomorphisms and isomorphisms,
relations between compactifications and completeness are scrutinized
over topological quasigroups. There are investigated applications to
families of maps over  topological quasigroups, restrictions of
mappings, their openness, continuous bijections. Necessary theorems
about topological quasigroups are proved. \footnote{key words and
phrases: regularity; completeness;
topological space; topological quasigroup; continuous mapping   \\
Mathematics Subject Classification 2010: 20N05; 22A30; 54C05; 54E15;
54H11}

\end{abstract}
\par Address: Dep. Appl. Mathematics, Moscow State Techn. Univ. MIREA,
\par av. Vernadsky 78, Moscow 119454, Russia; e-mail: sludkowski@mail.ru

\section{Introduction}

\par Studies of embeddings of topological spaces,
completeness, regularity, real complete spaces play great role in
topology and its applications \cite{adbtb,eng}. An important role
have properties of families of functions from a topological space in
the field of real numbers. They are used also for investigations of
topological vector spaces of real valued functions
\cite{fjor10,mntpscfb,rtvput21}. They are applied in functional
analysis, abstract harmonic analysis, theory of random functions,
mathematical physics \cite{nai,hew,reedsim}. On the other hand,
applications of nonassociative algebras near to quasigroups are
being developed recently in noncommutative analysis, noncommutative
geometry, operator theory, PDEs, elementary particle physics,
quantum field theory (see
\cite{guertzeb,lawmich,culbjglta07,eaksath90,ludqimqgsmj23,ludseabmmiop15,
ludhpdecv13,ludspdecv16,movsabqb} and references therein). According
to studies in the 20-th century an existence of a nontrivial
noncommutative geometry is equivalent to that of a corresponding
quasigroup (see \cite{movsabqb,sabininb99} and references therein).
\par Recently nonassociative algebras near to quasigroups were
utilized in investigations of slave boson decompositions in
superconductors \cite{dnasbdjglta07} in nonassociative quantum
mechanics \cite{dtmnaqmahep07}. They were also actively used in
gauge theories and Green-Schwartz superstrings
\cite{hdtjmp91,mgtnasjmp05}. Nonassociative algebras of such type
are connected with quasi-hopf deformations in nonassociative quantum
mechanics \cite{mssjmp14}. Nonassociative algebras near to
quasigroups served as one of the main tools during researches of De
Sitter representations of a curved space-time \cite{knasqmcstjmp99},
in the great unification theory and for studies of Yang-Mills fields
\cite{guertzeb,cgostjmp07}. The family of such nonassociative
algebras was utilized for an analysis of Yang-Baxter PDEs with
applications to the great unification theory (see
\cite{incsybea21,nichitaax19,ninybeut14} and references therein).
Quasigroups have found other applications in informatics and coding
theory, because they open new opportunities in comparison with
groups \cite{blautrctb,gkmngcnagdm04,mmnnascdjms20,pualdb94}. It is
necessary to emphasize that themes indicated at the beginning were
intensively investigated over the real field $\mathbf{R}$ and the
complex field $\mathbf{C}$, but very little is known about families
of functions with values in topological quasigroups. According to
the discussion presented above their study is stimulated, for
example, by functional analysis, abstract harmonic analysis, random
functions and stochastic processes, mathematical physics over
algebras near to quasigroups, by investigations of microbundle
structure \cite{ludkmtrta19}. The aforementioned domains are based
on algebro-topological binary systems.
\par Thus a new perspective direction for investigations opens
of topological quasigroups and their actions on topological spaces,
about which it is known very little in comparison with topological
groups and transformations of topological spaces by groups.
\par This article is devoted to completely regular and complete topological
spaces over topological quasigroups. Separations and embeddings are
also investigated for them. Their homeomorphisms and isomorphisms,
relations between compactifications and completeness are scrutinized
over topological quasigroups. There are investigated applications to
families of mappings over  topological quasigroups, restrictions of
maps, their openness, continuous bijections. Necessary theorems
about topological quasigroups are proved in subsection 1.1.

\subsection{Some specific features of topological quasigroups}
\par {\bf Theorem 1.1.} {\it Let $G=Q(\mathcal{A})$ be a topological quasigroup, where $Q$ is a
topological space, $\mathcal{A}: Q^2\to Q$ is a continuous binary
operation on $G$. Then $Q(\mathcal{A})$ is principally isotopic to a
topological unital quasigroup $G_1=Q(\mathcal{B})$.} \par {\bf
Proof.} Let $a$ and $b$ belong to $Q$. Put
$\mathcal{B}(x,y)=\mathcal{A}(R_a^{-1}x,L_b^{-1}y)$ for each $x$ and
$y$ from $Q$, where $R_ax=\mathcal{A}(x,a)$,
$L_bx=\mathcal{A}(b,x)$, $R_a^{-1}d=d/a$ is a unique solution of the
equation $\mathcal{A}(a,y)=d$, $L_b^{-1}d=b\setminus d$ is a unique
solution of the equation $\mathcal{A}(b,y)=d$, where $a$ and $b$ are
given, while $y$ is to be calculated. Since $G$ is the topological
quasigroup, then the mappings $Q^2\ni (a,d)\mapsto d/a\in Q$ and
$Q^2\ni (b,d)\mapsto b\setminus d\in Q$ are continuous. By virtue of
theorem 1.3.2 in \cite{movsabqb} $G_1=Q(\mathcal{B})$ is a unital
quasigroup with the unit element $e=\mathcal{A}(b,a)$,
$\mathcal{B}=\mathcal{A}^{(R_a^{-1},L_b^{-1}, id_Q )}$,
$Q(\mathcal{B})$ is a $LR$-isotop of the quasigroup
$Q(\mathcal{A})$, where $id_q(q)=q$ for each $q\in Q$. Notice that
the mappings $\mathcal{A}$, $R_d$, $R_d^{-1}$, $L_d$, $L_d^{-1}$ are
continuous for each $d\in G$, consequently, $\mathcal{B}: Q^2\to Q$
is continuous and $Q(\mathcal{B})$ is the topological unital
quasigroup.

\par {\bf Theorem 1.2.} {\it Let $G=Q(\mathcal{A})$ be a topological unital quasigroup,
${\sf C}_e$ be a connected component of the unit element $e$ in $G$,
then ${\sf C}_e$ is a closed invariant unital subquasigroup in $G$,
${\sf C}_eb=b{\sf C}_e$ is a connected component of $b\in G$.}
\par {\bf Proof}. Since $e\in {\sf C}_e$ and ${\sf C}_e$
is a maximal connected subset in $Q$ containing $e$, then $b{\sf
C}_e$ and ${\sf C}_eb$ are connected and contain $b\in G$, since
$L_b: G\to G$ and $R_b: G\to G$ are the homeomorphisms. Therefore
${\sf C}_eb=b{\sf C}_e$ is the connected component ${\sf C}_b$ of
$b\in G$, since $L_b^{-1}: G\to G$ and $R_b^{-1}: G\to G$ are the
homeomorphisms. Then ${\sf C}_{ab}=(ab){\sf C}_e$, $(ab){\sf
C}_e=a(b{\sf C}_e)$, $a(b{\sf C}_e)=a({\sf C}_eb)$, $a({\sf
C}_eb)=(a{\sf C}_e)b$, $(a{\sf C}_e)b=({\sf C}_ea)b$, $({\sf
C}_ea)b={\sf C}_e(ab)$, ${\sf C}_e(ab)={\sf C}_{ab}$. Notice that
$\forall a\in G, \forall b\in G, ({\sf C}_a\cap {\sf C}_b \ne
\emptyset )$ $\leftrightarrow $ $({\sf C}_a= {\sf C}_b)$, the
connected component ${\sf C}_b$ is closed in $G$ for each $b\in G$.
Therefore ${\sf C}_e{\sf C}_e$, ${\sf C}_e\setminus {\sf C}_e$,
${\sf C}_e/{\sf C}_e$ are connected and contain $e$, since ${\sf
C}_b/b={\sf C}_e$ and $b\setminus {\sf C}_b={\sf C}_e$ for each
$b\in G$, ${\sf C}_b={\sf C}_e$ for each $b\in {\sf C}_e$, where
$ab:=A(a,b)$ is a shortening of the notation corresponding to
multiplication for all $a\in G$ and $b\in G$. Thus ${\sf C}_e{\sf
C}_e={\sf C}_e$, ${\sf C}_e\setminus {\sf C}_e={\sf C}_e$, ${\sf
C}_e/{\sf C}_e={\sf C}_e$.

\par {\bf Theorem 1.3.} {\it Let $G=Q(\mathcal{A})$ be a topological unital quasigroup, ${\sf C}_e$
be a connected component of the unit element $e=e_G$. Then there
exists a totally disconnected $P=G/_c{\sf C}_e$  quotient $T_1\cap
T_3$ unital quasigroup.}
\par {\bf Proof.} By virtue of theorem 1.2 ${\sf C}_e$ is the
invariant closed unital subquasigroup in $G$. From corollary 2.2 in
\cite{ludumj23} it follows, that a space $P=G/_c{\sf C}_e$ of all
left classes $\{ b {\sf C}_e: b\in G \} $ is a topological unital
quasigroup and the quotient mapping $\pi : G\to P$ is a continuous
open homomorphism. Notice that $P$ is a $T_1\cap T_3$ unital
quasigroup according to theorem 2.3 in \cite{ludumj23}, since ${\sf
C}_e$ is closed in $G$. Consider a subset $X\subset G$ such, that
${\sf C}_e\subset \bigcup \{ b {\sf C}_e: b\in X \} =:Y$ and $Y\ne
{\sf C}_e$. If $J\subset G$, then $\pi (J\cap (X{\sf C}_e))=\pi
(J)\cap (X{\sf C}_e)$. Since $Y\ne {\sf C}_e$, then $Y$ is
disconnected, consequently, there exist nonvoid open subsets $U_1$
and $U_2$ in $G$ such, that $Y=(U_1\cap Y)\cup (U_2\cap Y)$ and
$(U_1\cap Y)\cap (U_2\cap Y)=\emptyset $. Therefore $X{\sf C}_e=
(\pi (U_1)\cap (X{\sf C}_e))\cup (\pi (U_2)\cap (X{\sf C}_e))$,
where $\pi (U_1)$ and $\pi (U_2)$ are open in $Q$, since $\pi $ is
the continuous open homomorphism. For each $b\in G$ the left class
$b{\sf C}_e={\sf C}_b$ is connected, consequently, either $b{\sf
C}_e\subset U_1\cap Y$ or $b{\sf C}_e\subset U_2\cap Y$. Then $(\pi
(U_1)\cap (X{\sf C}_e))\cap (\pi (U_2)\cap (X{\sf C}_e))=\emptyset
$, consequently, $\pi (X)=X{\sf C}_e$ is disconnected.  Thus the
connected component of the unit element $e_P$ in $P$ coincides with
$e_P$.

\par {\bf Theorem 1.4.} {\it Let $G=Q(\mathcal{A})$ be a topological $T_1$ quasigroup.
Then it possesses a uniformity $\mathcal{U}_Q$ compatible with its
topology.}
\par {\bf Proof.} By virtue of theorem 1.1 it is sufficient to consider
a topological $T_1$ unital quasigroup $G_1=Q(\mathcal{B})$
principally isotopic to the quasigroup $G$. Let $\mathcal{T}_Q$ be a
topology on $Q$. Then the mappings $L_g: Q\to Q$ and $L_g^{-1}: Q\to
Q$ are homeomorphisms for each $g\in G_1$. \par We consider a
minimal group $H$ generated by homeomorphisms $L^{\epsilon
_1}_{g_1}\circ ... \circ L^{\epsilon _k}_{g_k}: (Q, \mathcal{T}_Q)
\to (Q, \mathcal{T}_Q)$ for all possible $k\in \mathbf{N}$, $g_j\in
G_1$, $\epsilon _j \in \{ -1, 1 \} $, $j\in \{ 1,..., k \} $, where
$L^{\epsilon _1}_{g_1}\circ L^{\epsilon _2}_{g_2}(g)=L^{\epsilon
_1}_{g_1}(L^{\epsilon _2}_{g_2}(g))$, $L_{g_2}(g)=B(g_2,g)$,
$L_{g_2}^{-1}(g)$ is a unique solution of the equation $B(g_2,y)=g$
with $y$ to be calculated and for given $g$ and $g_2$ in $G_1$.
Supply $H$ by a base of a topology $\mathcal{B}_H = \{
W(g_1,...,g_k;L_{g_{k+1}}(V)): k\in \mathbf{N}; g_j\in G_1; j\in \{
1,..., k+1 \}; V \in \mathcal{B}_{e,G_1} \} $, where $W(g_1,...,g_k;
L_{g_{k+1}}(V)) := \{ g\in G_1: L_{g_j}^{\epsilon _j}(g)\in
L_{g_{k+1}}(V); \epsilon _j \in \{ -1, 1 \}; j\in \{ 1,..., k \} \}
$,  $\mathcal{B}_{e,G_1}$ is a base of open neighborhoods of the
unit element $e$ in $G_1$. Since $G_1$ is the topological unital
quasigroup, then for each open neighborhood $V$ of the unit element
$e$ in $G_1$, $n\in \mathbf{N}$, $\epsilon _j \in \{ -1, 1 \} $,
$j\in \{ 1,..., n \} $, there exists an open neighborhood $W$ for
$e$ such that $L^{\epsilon _1}_{W_1}\circ ... \circ L^{\epsilon
_k}_{W_k}\subset V$, where $L_W^{\epsilon }= \{ L_g^{\epsilon }:
g\in W \} $. Therefore $(H, \mathcal{T}_H)$ is the topological
group, where $\mathcal{T}_H$ is the topology on $H$ with the base
$\mathcal{B}_H$.
\par A mapping $\phi : (Q,\mathcal{T}_Q) \to (H,\mathcal{T}_H)$ such that $\phi (q)=L_q$
for each $q\in Q$ is continuous by the construction given above, it
separates points and closed subsets in $Q$, since $G_1$ is the
unital quasigroup, $B(q,P)$ is closed in $(Q,\mathcal{T}_Q)$ for
each closed $P$ in $(Q,\mathcal{T}_Q)$. From theorem 2.3.20 and
lemma 2.3.19 in \cite{eng} it follows, that $\phi :
(Q,\mathcal{T}_Q) \to (\phi (Q), \mathcal{T}_H|_Q) \subseteq
(H,\mathcal{T}_H)$ is the homeomorphism. On the other hand, on the
topological group $(H,\mathcal{T}_H)$ a left uniformity exists
$\mathcal{U}_H$ compatible with the topology $\mathcal{T}_H$ on it
(see (4.14)(b), definitions in \cite{hew} and example 8.1.17 in
\cite{eng}). Then $\mathcal{U}_H$ induces a uniformity
$\mathcal{U}_Q$ on $Q$ compatible with the topology $\mathcal{T}_Q$
on $Q$.

\par {\bf Remark 1.1.}Apart from a topological group a uniformity
on a quasigroup generally may be neither left nor right invariant,
since the mapping $\phi $ may be not a homomorphism, since $G$ may
be nonassociative being a quasigroup.

\section{Regularity and completeness over topological quasigroups}
\par {\bf Definition 2.1.} Let $X$ be a topological space,
$G=Q(\mathcal{A})$ be a topological nontrivial infinite quasigroup.
Let $G$ be a $T_1$ topological space with a topology
$\mathcal{T}_G$, and $G$ be complete relative to its own uniformity
$\mathcal{U}_G$. If $X$ is a $T_1$ space, and for each closed subset
$B$ and each point $x$ in $X$ such that $x\notin B$, there exists
$q\ne r$ in $Q$ and a continuous mapping $f: X\to Q$ such that
$f(x)=r$ and $f(B)= \{ q \} $, then $X$ will be called a weakly
completely $G$-regular space. If it is accomplished for fixed $q\ne
r$ in $Q$ and every $x$ and $B$ satisfying the conditions given
above, then $X$ will be called a completely $(G,q,r)$-regular space.
If $X$ completely $(G,q,r)$-regular for each $q\ne r$ in $Q$, then
$X$ will be called a completely $G$-regular space.

\par {\bf Remark 2.1.} Definition 2.1 means that
if the topological space $X$ is completely $(G,q,r)$-regular, then
it is weakly completely $G$-regular; if $X$ is completely
$G$-regular, then $X$ is completely $(G,q,r)$-regular. In
particular, if $X$ is a Tychonoff zero-dimensional space, then $X$
is completely $G$-regular. Indeed, for $x$ there exists a clopen
(that is open-and-closed) neighborhood $V\ni x$ with $V\subset
X\setminus B$. As $f$ it is possible to choose the characteristic
function $\chi _{X\setminus V}$, where $\chi _A(y)=q$ for each $y\in
A$, $\chi _A(y)=r$ for each $y\notin A$, for $A\subset X$, since
$q\ne r$, $G=Q(\mathcal{A})$ is the topological $T_1$ quasigroup.
The mapping $\chi _{X\setminus V}$ is continuous, since $V$ is
clopen in $X$.
\par By virtue of theorem 2.1 given below the topological quasigroup
$G$ itself is weakly completely $G$-regular.

\par {\bf Theorem 2.1.}
{\it Let $G=Q(\mathcal{A})$ be a topological $T_1$ quasigroup
complete relative to its uniformity. Then $G$ is weakly completely
$G$-regular.}
\par {\bf Proof.} By virtue of theorem 1.1 it sufficient to consider
the case of the topological unital quasigroup $G_1=Q(\mathcal{B})$
with $r=e$. By theorem 1.2 the connected component ${\sf C}_e$ of
the unit element is the closed  invariant unital subquasigroup in
$G$. Moreover, ${\sf C}_eg=g{\sf C}_e$ is the connected component of
an arbitrary element $g\in G_1$. Notice that the quotient unital
quasigroup $P={G_1}/_c{\sf C}_e$ is $T_2$ and totally disconnected
according to theorem 1.3.
\par If ${\sf C}_e$ is nontrivial, then for the topological
space $g{\sf C}_e$ there exists a decomposition into the limit of an
inverse spectrum of polyhedra over the real field ${\mathbf R}$ by
Freudenthal theorem for each $g\in G_1$ (see reference in \cite{eng}
or in \cite{ludnpr2000}). If the quotient unital quasigroup
$P=G_1/_c{\sf C}_e$ is nontrivial, then the topological space
$Q/_c{\sf C}_e$ has the decomposition into the limit of the inverse
spectrum of polyhedra over the field $\mathbf{L}\supseteq
\mathbf{Q}_p$, which is a finite algebraic extension of the field
$\mathbf{Q}_p$ of $p$-adic numbers
 (see theorem 3.19 in \cite{ludnpr2000}). The polyhedra over
${\mathbf R}$ are $T_1\cap T_{3\frac{1}{2}}$ spaces, the polyhedra
over $\mathbf{L}$ are totally disconnected.
\par For a closed subset $B$ and a point $x$ in $Q$ such that
$x\notin B$, we choose a point $q$ in $B$. For each $b\in Q$ the
left shift mapping $L_b$ is the homeomorphism from $Q$ onto $Q$ as
the topological space, where $L_by=by$ for each $b$ and $y$ in $Q$.
Therefore without restriction of generality it is possible to choose
$x=e$. From the polyhedral decompositions it follows that there
exists a continuous mapping $f: Q\to Q$ such that $f(x)=e$, $f(B) =
\{ q \} $.

\par {\bf Corollary 2.1.} {\it If $G$ is the topological $T_1$ quasigroup complete
relative to its uniformity either connected or totally disconnected,
then it is completely $G$-regular.}

\par {\bf Remark 2.2.} Besides a quasigroup it is possible to take an algebra
$\mathbf{A} = Q(\mathcal{A}_1, \mathcal{A}_2)$ with two binary
operations $\mathcal{A}_1$ and $\mathcal{A}_2$ from $Q^2$ onto $Q$
such that $Q(\mathcal{A}_1)$ is a topological quasigroup and $(Q- \{
x_0 \}) (\mathcal{A}_2)$ is a topological right quasigroup,
$Q(\mathcal{A}_2)$ is a topological groupoid, $\mathcal{A}_2(Q,x_0)=
\{ x_0 \} $, where $x_0$ is a marked point in $Q$.

\par {\bf Corollary 2.2.} {\it If $\mathbf{A}$ is the topological $T_1$ algebra
as in remark 2.2, then it is completely $\mathbf{A}$-regular.}
\par {\bf Proof.} By virtue of theorem 2.1 $\mathbf{A}$ is weakly completely
$G$-regular, where $G=Q(\mathcal{A}_1)$. Using the left shift
operator $L_b$ on $G$ without loss of generality it is possible to
choose $r=x_0$ and $q\ne r$, where $L_bx=\mathcal{A}_1(b,x)$ for
each $x\in G$. Since $(Q- \{ x_0 \}) (\mathcal{A}_2)$ is the right
quasigroup, then for each $a$ and $b$ in $(Q- \{ x_0
\})(\mathcal{A}_2)$ the equation $\mathcal{A}_2(y,a)=b$ has a unique
solution $y$ in $(Q- \{ x_0 \})(\mathcal{A}_2)$. Then for each $a\ne
x_0$ in $Q$ the right shift mapping $R_a: Q(\mathcal{A}_2) \to
Q(\mathcal{A}_2)$ is continuous, where $R_ax=\mathcal{A}_2(x,a)$ for
each $x$ in $Q$. Notice that if the mapping $f: Q\to Q$ is
continuous, then the mapping $\mathcal{A}_2(y,f): Q\to Q$ is
continuous, where $\mathcal{A}_2(y,f)(x)=\mathcal{A}_2(y,f(x))$ for
each $x\in \mathbf{A}$, $y\in \mathbf{A}$.

\par {\bf Corollary 2.3.} {\it Let $G$ be the topological $T_1$ quasigroup
complete relative to its uniformity either connected or totally
disconnected; or $G=Q(\mathcal{A}_1, \mathcal{A}_2)$ be the
topological $T_1$ algebra satisfying the conditions of remark 2.2.
If a space $X$ is weakly completely $G$-regular, then it is
completely $G$-regular.}

\par {\bf Theorem 2.2.} {\it Let $A$ be a subspace of a (weakly) completely $G$-regular
or $(G,q,r)$-regular space $X$. Then $A$ is (weakly) completely
$G$-regular or $(G,q,r)$-regular correspondingly.}
\par {\bf Proof.} By virtue of proposition 2.1.1 in \cite{eng} and definition
2.1 $A$ is the $T_1$ space. Let $x\in A$, $B$ be a closed subset in
$A$, $x\notin B$. By proposition 2.1.1 in \cite{eng} $B=A\cap
\bar{B}$, consequently, $x\notin \bar{B}$, where $\bar{B}=cl_XB$
denotes the closure of $B$ in $X$. In the case of the weak complete
$G$-regularity there exist $q\ne r$ in $Q$ and a continuous mapping
$f: X\to Q$ such that $f(x)=r$ and $f(\bar{B})=\{ q \} $. If $X$ is
completely $(G,q,r)$-regular, then a continuous mapping $f: X\to Q$
exists such that $f(x)=r$ and $f(\bar{B})=\{ q \} $. In the case of
the complete $G$-regularity for each $q\ne r$ in $Q$ a continuous
mapping $f: X\to Q$ exists such that $f(x)=r$ and $f(\bar{B})=\{ q
\} $. Then the restriction of $f$ on $A$ is continuous $f|_A: A\to
Q$, with $f(x)=r$ and $f(B)= \{ q \} $, consequently, $A$ is weakly
completely $G$-regular or completely $G$-regular or
$(G,q,r)$-regular correspondingly.

\par {\bf Proposition 2.1.} {\it Let $X$ be a $T_1$ space,
$\mathcal{P}_X$ be a subbase of its topology $\mathcal{T}_X$. $X$ is
completely $G$-regular if and only if for each $x\in X$ and each
neighborhood $V$ of a point $x$, where $V$ belongs to the subbase
$\mathcal{P}_X$, for each $q\ne r$ in $G$ a continuous mapping $f:
X\to G$ exists such that $f(x)=r$ and $f(X\setminus V)=\{ q \} $.}
\par {\bf Proof.} The necessity follows from definition 2.1,
since $x\notin X\setminus V$, $X\setminus V$ is closed in $X$.
\par Sufficiency. Let $x\in X$, $B$ be closed in $X$, $x\notin B$.
By the definition of the subbase there exist $V_1,...,V_k$ in
$\mathcal{P}_X$ such that $x\in \bigcap_{i=1}^k V_i$ and
$\bigcap_{i=1}^k V_i\subseteq X\setminus B.$ According to definition
2.1 for each $i\in \{ 1,...,k \} $ and for each $q_i\ne r_i$ in $G$
a continuous mapping $f_i: X\to G$ exists such that $f_i(x)=r_i$ and
$f_i(X\setminus V_i)= \{ q_i \} $. Put $f(z)=f_1(z)$ for $k=1$,
$f(z)=f_1(z)f_2(z)$ for $k=2$, $f(z)=(...(f_1(z)f_2(z))...)f_k(z)$
for $k\ge 3$, for each $z\in X$. Then $f$ is the continuous mapping,
since $G$ is the topological quasigroup. Notice that $f(x)=r$ and
$f(y)=q$ for each $y\in B$, since $B\subseteq X\setminus V_i$ for
each $i=1,...,k$, where $r=r_1$ and $q=q_1$ for $k=1$, $r=r_1r_2$
and $q=q_1q_2$ for $k=2$, $r=(...(r_1r_2)...)r_k$ and
$q=(...(q_1q_2)...)q_k$ for $k\ge 3$. Since $G$ is the quasigroup,
then for each $q\ne r$ in $G$ there exist $q_1\ne r_1$,..., $q_k\ne
r_k$ in $G$ such that $q$ and $r$ have decompositions into the
ordered products provided above.

\par {\bf Theorem 2.3.} {\it Let $X=\prod_{j\in J}X_j$ be the product of topological
spaces, where $J$ is a set, $X\ne \emptyset $. $X$ is completely
$G$-regular if and only if $X_j$ is completely $G$-regular for each
$j\in J$.}
\par {\bf Proof.} Let $X_j$ be completely $G$-regular
for each $j\in J$. By virtue of proposition 2.1 it is sufficient to
prove, that for each point $x\in X$ and each its neighborhood of the
form $V=\pi ^{-1}_{j_0}(W_{j_0})$ and for each $q\ne r$ in $G$,
where $W_{j_0}$ is open in $X_{j_0}$, $j_0\in J$, $\pi _j: X\to X_j$
is a projection, there exists a continuous mapping $f: X\to G$ such
that $f(x)=r$ and $f(X\setminus V)= \{ q \} $, since $\prod_{j\in
J}X_j$ is supplied with the Tychonoff product topology  (see also
section 2.3 in \cite{eng}). From the complete $G$-regularity of the
space $X_{j_0}$ according to definition 2.1 it follows, that there
exists a continuous mapping $g: X_{j_0}\to G$ such that
$g(x_{j_0})=r$ and $g(X_{j_0}\setminus W_{j_0})= \{ q \} $, where
$x_{j_0}=\pi _{j_0}(x)$. Then the mapping $f(y)=g(\pi _{j_0}(y))$
satisfies the demanded conditions, where $y\in X$.
\par Vice versa, let $X$ be completely $G$-regular. Choose
an arbitrary fixed point $y_j\in X_j$ for each $j\in J$. Consider
the product topological space $Y_{j_0}=\prod_{j\in J}A_j$, where
$A_{j_0}=X_{j_0}$, $A_j= \{ y_j \} $ for each $j\in J\setminus \{
j_0 \} $. Since $\pi _{j_0}: Y_{j_0}\to X_{j_0}$ is the
homeomorphism, then there exists an embedding of $X_{j_0}$ into
$Y_{j_0}$. By definition 2.1 $X_j$ is the $T_1$ space for each $j\in
J$. Therefore $Y_{j_0}$ is closed in $X$. Thus, there exists an
embedding of $X_{j_0}$ into $X$ as the closed subspace,
consequently, $X_{j_0}$ is completely $G$-regular by virtue of
theorem 2.2.

\par {\bf Theorem 2.4.} {\it Let $X$ be weakly completely $G$-regular space,
$w(X)={\sf m}\ge \aleph _0$, where $w(X)$ denotes weight of the
topological space $X$. Then there exists a homeomorphic embedding of
$X$ into $G^{\sf m}$.}
\par {\bf Proof.} From theorem 2.3.13 in \cite{eng} it follows,
that $w(G^{\sf m})\le {\sf m}w(G)$, since ${\sf m}w(G)=\max ({\sf
m}, w(G))$, ${\sf m}{\sf m}={\sf m}$, $w(G)w(G)=w(G)$, ${\sf m}\ge
\aleph _0$, $w(G)\ge \aleph _0$ \cite{kunenb}. There exists an
embedding of the discrete space $D({\sf m})$ of the cardinality
${\sf m}$ into $G^{\sf m}$, since $G$ is the infinite $T_1$
topological quasigroup. From $w(D({\sf m}))={\sf m}$ it follows,
that $w(G^{\sf m})\ge {\sf m}$. By virtue of theorem 1.1.15 in
\cite{eng} there exists a base $ \{ U_j: j\in J \} $ in the space
$X$ consisting of functionally open sets $U_j=f_j^{-1}(V_j)$, where
$f_j: X\to G$ is a continuous mapping for each $j\in J$, $V_j$ is an
open subset in $G$, $J$ is a set of the cardinality ${\sf m}$.
\par If $x\in X$, $B$ is closed in $X$, $x\notin B$, then
$W=X\setminus B$ is open in $X$, $x\in W$, then there exists $j_0\in
J$ such that $U_{j_0}\subseteq W$, consequently, $f_{j_0}$ separates
$x$ and $B$. Thus, the family ${\mathcal F}=\{ f_j: j\in J \} $
separates points and closed sets in $X$. Since $X$ is the $T_1$
space, then by theorem about diagonal mapping 2.3.20 in \cite{eng}
the diagonal $\Delta _{j\in J}f_j$ is the homeomorphic embedding of
$X$ into $G^{\sf m}$.

\par {\bf Remark 2.2.} If in the definition of the complete regularity
use instead of the entire topological quasigroup $G$ its canonical
closed subset $W$ in $G$, then there will be the definition of a
(weakly) completely $W$-regular ($(W,q,r)$-regular) space. Then
theorems 2.2 and 2.4 can be transferred for (weakly) completely
$W$-regular spaces with $f: X\to W$ and $W^{\sf m}$ instead of $f:
X\to G$ and $G^{\sf m}$ correspondingly. If additionally
$WW\subseteq W$, then proposition 2.1 and theorem 2.3 also can be
transferred for completely $W$-regular spaces.

\par {\bf Definition 2.2.} A topological space is called
$G$-complete, if it is weakly completely $G$-regular and there does
not exist a weakly completely $G$-regular space $Y=Y_X$, which
satisfies the following conditions: \par $(i)$ there exists a
homeomorphic embedding $h: X\to Y$ such that $h(X)\ne \overline{
h(X)} =Y$; \par $(ii)$ for each continuous mapping $f: X\to G$ there
exists a continuous mapping $g: Y\to G$ such that $g(h)=f$, where
$\bar{A}=cl_YA$ denotes the closure of a subset $A$ in $ Y$.

\par {\bf Theorem 2.5.} {\it A topological space $X$ is $G$-complete
if and only if $X$ is homeomorphic to a closed subspace in $G^J$ for
some set $J$, where $G^J$ is supplied the Tychonoff product
topology.}
\par {\bf Proof.} Let $X$ be closed in $G^J$ and
a homeomorphic embedding $h: X\to Y$ exist, where $Y$ is a weakly
completely $G$-regular space. Let $h$ satisfy condition $(ii)$ of
definition 2.2. Without loss of generality it is possible to
consider the case $cl_Yh(X)=Y$, where $cl_YA$ denotes the closure of
a subset $A$ in $Y$. By virtue of condition $(ii)$ for each
projection $\pi _j: G^J\to G_j$, where $G_j=G$, $G^J=\prod_{j\in J}
G_j$, there exists a continuous mapping $\eta _j: Y\to G_j$ such
that $\eta _j(h)=\pi _j|_X$, since each projection $\pi _j$ has the
continuous restriction $\pi _j|_X$ on $X$. Therefore there exists a
continuous injective diagonal mapping $\eta =\Delta _{j\in J}\eta
_j: Y\to G^J$ by virtue of theorem 2.3.20 in \cite{eng}. Then $\eta
(h(x))=x$ for each $x\in X$, consequently, $\eta (Y)=\eta
(cl_Yh(X))\subseteq cl_{G^J}\eta (h(X))=X$, since $X$ is closed in
$G^J$. Notice that $h(\eta (h(x)))=h(x)$ for each $x\in X$. Thus,
$h(\eta |_{h(X)})=id_{h(X)}$. Since the topological $T_1$ quasigroup
$(G, \mathcal{T}_G)$ possesses a uniformity $\mathcal{U}_G$
compatible with its topology according to theorem 1.4, then it is
the topological $T_1\cap T_{3\frac{1}{2}}$ space by theorem 8.1.20
in \cite{eng}. By virtue of theorem 2.3.11 in \cite{eng} $G^J$ is
the $T_1\cap T_{3\frac{1}{2}}$ space. From theorem 1.5.4 in
\cite{eng} it follows, that $h(\eta |_{h(X)})=id_Y$. Thus, $h(X)=Y$,
since $h(\eta (Y))=h(X)$. Consequently, there does not exist the
weakly completely $G$-regular space satisfying condition $(i)$ and
$(ii)$ in definition 2.2. This means that $X$ is $G$-complete.
\par Let now $X$ be the $G$-complete topological
space. Choose a space $C(X,G)$ of all continuous mappings $f: X\to
G$ and the diagonal mapping $\psi := \Delta _{f\in C(X,G)}f: X\to
\prod_{f\in C(X,G)}G_f$, where $G_f=G$ for each $f\in C(X,G)$. The
family $C(X,G)$ separates points and closed sets in $X$ according to
definitions 2.1 and 2.2. By virtue of theorem 2.3.20 in \cite{eng}
about the diagonal mapping $\psi $ is the homeomorphic embedding.
Then take $Y=cl_{G^J}\psi (X)$, where $J=C(X,G)$. For each
continuous mapping $g: X\to G$ there exists a continuous mapping
$\tilde{g}: Y\to G$ such that $\tilde{g}(\psi )=g$ and
$\tilde{g}=\pi _g|_Y$, where $\pi _g: \prod_{f\in J}G_f\to G_g$.
From definition 2.2 it follows, that $Y=\psi (X)$. Thus, $X$ is
homeomorphic to the closed subspace $\psi (X)$ in $G^J$.

\par {\bf Lemma 2.1.} {\it Let $G=Q(\mathcal{A})$ be an infinite topological nondiscrete
$T_1$ quasigroup. Let there exist continuous mappings $h_j: G\to G$
such that $\mathcal{A}(h_1,h_2)=id$; $h_j(G)\ne G$, $Int_Gh_j(G)$ is
$G$-complete for each $j\in \{ 1, 2 \} $, $Int_G(h_1(G)\cap
h_2(G))=\emptyset ,$ $cl_G(Int_Gh_1(G)\cup Int_Gh_2(G))=G$, where
$\forall x\in G$ $id(x)=x$. Let also $X$ be a topological space,
$A\subset X$, each continuous mapping $f: A\to G$ with $f(A)\subset
cl_G(Int_GW_j)$ for some $j\in \{ 1, 2 \} $ has a continuous
extension $F: G\to G$, where $W_1=h_1(G)$, $W_2=h_2(G)$. Then each
continuous mapping $f: A\to G$ has a continuous extension on $X$. }
\par {\bf Proof.} Since the mappings $h_1$ and $h_2$
are continuous, then the mapping $\mathcal{A}(h_1,h_2)$ is also
continuous on $G$, since
$\mathcal{A}(h_1,h_2)(x)=\mathcal{A}(h_1(x),h_2(x))$ for each $x\in
G$. From the conditions of this lemma it follows, that $G$ is dense
in itself and $cl_G(Int_GW_1)\cup cl_G(Int_GW_2)=G.$ Put
$f_1(x)=h_1(f(x))$, $f_2(x)=h_2(f(x))$ for each $x\in A$. Then $f_j:
A\to G$ and $f_j(A)\subseteq cl_G(Int_GW_j).$ By the condition of
this lemma there exists the continuous extension $F_j$ on $G$ for
$f_j$, $F_j: G\to G$, $F_j|_A=f_j$, for each $j\in \{ 1, 2 \} $.
Therefore $F=\mathcal{A}(F_1,F_2)$ is the continuous extension on
$G$ of the mapping $f$.

\par {\bf Example 2.1.} If $G=(\mathbf{R},+)$, then $h_1(x)=\max (x,0)$,
$h_2(x)= \min (x,0)$, $h_1(x)+h_2(x)=x$ for each $x\in \mathbf{ R}$,
$h_1(\mathbf{R})=[0, + \infty )$, $h_2(\mathbf{ R})=(-\infty ,0]$.
If a topological field $\mathbf{F}$ is an extension of the real
field $\mathbf{R}$, then $\mathbf{F}$ has the structure of the
topological vector space over $\mathbf{R}$. Therefore in this case
also there exist functions $h_1$ and $h_2$, satisfying the
conditions of lemma 2.1.

\par {\bf  Example 2.2.} Let $G=Q(\mathcal{B})$ be a totally disconnected
topological infinite nondiscrete $T_1$ unital quasigroup. For any
clopen nonvoid subsets $A_1$ and $A_2$ in $G$ such that $A_1\cup
A_2=G$, $A_1\cap A_2=\emptyset $, it is possible to put $h_1(x)=x$
for $x\in A_1$, $h_1(x)=e$ for $x\in A_2$, and $h_2(x)=x$ for $x\in
A_2$, $h_2(x)=e$ for $x\in A_1$, where $e$ is the unit element in
$G$. Then mappings $h_1$ and $h_2$ are continuous,
$\mathcal{B}(h_1(x),h_2(x))=x$ for each $x\in G$, $h_1(G)=A_1\cup \{
e \} $, $h_2(G)=A_2\cup \{ e \} $. \par Using theorem 1.1 the
mappings $h_1$ and $h_2$ can be constructed on the totally
disconnected topological infinite nondiscrete $T_1$ quasigroup
$Q(\mathcal{A})$.

\par {\bf Theorem 2.6.} {\it Let
$X$ be a weakly completely $G$-regular space. Then the following
conditions are equivalent: \par $(i)$ $X$ is $G$-complete;
\par $(ii)$ for each $j\in \{ 1,2 \} $ there exists a compactification
$Z_j=c(cl_G(Int_G(W_j)))$ such that for each $y\in (\beta
X)\setminus X$ there exists $i=i(y)\in \{ 1, 2 \} $, and there
exists a continuous mapping $h_{i,y}: \beta X\to Z_i$ such that
$h_{i,y}(y)\in Z_i\setminus ((cl_G(Int_G(W_i)))\setminus (W_1\cap
W_2))$, and $h_{i,y}(x)\in (cl_G(Int_G(W_i)))\setminus (W_1\cap
W_2)$ for each $x\in X$, where $W_1$ and $W_2$ satisfy the
conditions of lemma 2.1.}
\par {\bf Proof.} Let the topological space $X$ be $G$-complete, $y\in (\beta
X)\setminus X$. Put $Y=X\cup \{ y \} $, consequently, $Y\subseteq
\beta X$. Therefore the condition $(i)$ is satisfied of definition
2.2. From the condition $(ii)$ of definition 2.2 it follows, that
there exists a continuous mapping $f: X\to G$, which has not a
continuous extension on $Y$. By virtue of lemma 2.1 it is sufficient
to consider the case: there exists $j\in \{ 1, 2 \} $ such that
$f(X)\subseteq cl_G(Int_GW_j).$ By virtue of theorem 1.4 $G$
possesses the uniformity compatible with its topology, consequently,
it is $T_1\cap T_{3\frac{1}{2}}$ (that is Tychonoff) by theorem
8.1.20 in \cite{eng}. Since $G$ is $T_1\cap T_{3\frac{1}{2}}$, and
$X$ is weakly completely $G$-regular, then $X$ is $T_1\cap
T_{3\frac{1}{2}}$. From theorem 2.1 using the inverse spectrum
polyhedral decomposition in its proof it follows, that there exist
subsets $W_1$, and mappings $h_1$, $h_2$, satisfying the conditions
of lemma 2.1. By virtue of theorems 3.1.9 and 3.6.1 in \cite{eng}
the mapping $f$ has the continuous extension $\tilde{f}: \beta X\to
Z_j$. If there could be $f(y)\in cl_G(Int_GW_j)$, then there would
exist a continuous extension $\tilde{f}|_Y: Y\to G$ for $f$. It
provides the contradiction. Therefore $\tilde{f}(y)\in Z_j\setminus
cl_G(Int_GW_j)$.
\par Since $cl_G(V_1\cup V_2)=G$, where
$V_1=Int_Gh_1(G)$ and $V_2=Int_G(h_2(G))$, then $Z_1\cup Z_2\supset
G$. Let for each $y\in (\beta X)\setminus X$ there exist a
continuous mapping $h_{i,y}: \beta X\to Z_i$ such that
$h_{i,y}(y)\in Z_i\setminus ((cl_GInt_GW_i)\setminus (W_1\cap W_2))$
and $h_{i,y}(X)\subset (cl_GInt_GW_i)\setminus (W_1\cap W_2)$, where
$i=i(y)\in \{ 1, 2 \} $. Then $X=\bigcap _{y\in (\beta X)\setminus
X} h_{i(y),y}^{-1}((cl_GInt_GW_i)\setminus (W_1\cap W_2))$.
\par For the continuation of the proof there are necessary the following lemmas and
theorems and corollaries.

\par {\bf Lemma 2.2.} {\it Each closed subspace $X$ of the $G$-complete
space is $G$-complete.}
\par {\bf Proof.} Let $X$ be closed in the $G$-complete
space $Y$. By virtue of theorem 2.5 $Y$ is homeomorphic to a closed
subspace in $G^J$ for some set $J$. Therefore $X$ is homeomorphic to
a closed subspace in $G^J$, consequently, $X$ is $G$-complete by
theorem 2.5.

\par {\bf Theorem 2.7.} {\it Let $X_d\ne \emptyset $ for each $d\in D$, where
$D$ is a set. The product $\prod_{d\in D}X_d$ is $G$-complete if and
only if $X_d$ is $G$-complete for each $d\in D$.}
\par {\bf Proof.} If $X_d$ is $G$-complete, then $X_d$
is homeomorphic to a closed subspace in $G^{J_d}$ for some set $J_d$
by theorem 2.2. Therefore $\prod_{d\in D}X_d$ is homeomorphic to a
closed subspace in $G^J$ for $J=\bigcup_{d\in D}J_d$ by virtue of
proposition 2.3.7  and corollary 2.3.4 in \cite{eng}. From theorem
2.4 it follows, that $\prod_{d\in D}X_d$ is $G$-complete.
\par Vice versa, let $\prod_{d\in D}X_d$ be $G$-complete.
From corollary 2.3.4, proposition 2.3.7 in \cite{eng} and theorem
2.4 above it follows, that $X_d$ is homeomorphic to the closed
subspace in $G^{J_d}$ for some set $J_d$ for each $d\in D$,
consequently, $X_d$ is $G$-complete according to lemma 2.2.
\par {\bf Corollary 2.4.} {\it The limit of an inverse spectrum of
$G$-complete spaces is $G$-complete.} \par {\bf Proof.} It follows
from proposition 2.5.1 in \cite{eng}, lemma 2.2 and theorem 2.7
given above.
\par {\bf Theorem 2.8.} {\it Let $X$ be a topological space, $ \{ A_j: j\in J \} $ be
a family of its subspaces. If $A_j$ is $G$-complete for each $j\in
J$, then $\bigcap_{j\in J} A_j$ is $G$-complete.}
\par {\bf Proof.} The intersection $\bigcap_{j\in J}A_j$
is homeomorphic to $(\prod_{j\in J}A_j)\cap \Delta _X$, where
$\Delta _X$ is the diagonal in $\prod_{j\in J}X_j$, $X_j=X$ for each
$j\in J$. From theorem 2.7 it follows, that $\prod_{j\in J}A_j$ is
$G$-complete. By virtue of lemma 2.2 $\bigcap_{j\in J} A_j$ is
$G$-complete. \par {\bf Corollary 2.5.} {\it Let a topological space
$X$ be $G$-complete, $Y$ be a $T_2$ space, $f: X\to Y$ be a
continuous mapping, $D$ be as $G$-complete subspace in $Y$. Then the
inverse image $f^{-1}(D)$ is $G$-complete.}
\par {\bf Proof.} Consider the restriction $h$ of the mapping $f$ on
$f^{-1}(D) =\{ x\in X: f(x)\in D \} $ and $\Gamma (h)=(X\times
D)\cap \Gamma (f)$, where $\Gamma (f)= \{ (x,y)\in X\times Y: y=f(x)
\} $ is the graph of the mapping $f$. Then $\Gamma (h)$ is closed in
$X\times D$, and $X\times D$ is $G$-complete by theorem 2.7. Since
the inverse image $f^{-1}(D)$ is homeomorphic to $\Gamma (h)$, then
$f^{-1}(D)$ is $G$-complete according to lemma 2.2.
\par The continuation of the {\bf proof} of theorem 2.6. By virtue of
corollaries 2.4 and 2.5 $X$ is $G$-complete, since
$(cl_GInt_GW)\setminus (W_1\cap W_2)$ is $G$-complete, since
$Int_Gh_j(G)$ is $G$-complete for each $j\in \{ 1, 2 \} $.

\par {\bf Definition 2.3.} For a weakly completely $G$-regular
space $X$ let on the space $C(X,G)$ of all continuous mappings $f:
X\to G$ a subbase be given ${\sf B}(C(X,G)) = \{ W(x_1,...,x_n;U):
x_i\in X, i=1,...,n, n\in {\mathbf N}; U\in {\sf B}(G) \} $ of the
topology $\mathcal{ T}=\mathcal{ T}(C(X,G))$ on $C(X,G)$, where
$W(x_1,...,x_n;U) = \{ g\in C(X,G): g(x_i)\in U, i=1,...,n \} $,
$n\in \mathbf{ N}$, $\mathbf{ N}= \{ 1, 2, 3,... \} $, ${\sf B}(G)$
is the base of the topology of the quasigroup $G$. Then $(C(X,G),
\mathcal{ T})$ is called the space of all continuous mappings from
$X$ into $G$ in the topology of pointwise convergence and it is
denoted by $C_p(X,G)$.

\par {\bf Proposition 2.2.} {\it Let $X$ be a completely $G$-regular space. Then $C_p(X,G)$ is
an everywhere dense subspace in $G^X$.}
\par {\bf Proof.} From the complete $G$-regularity
of the space $X$ it follows, that for each $x_1,...,x_n$ in $X$ and
$f\in G^X$ there exists a mapping $h\in C_p(X,G)$ such that
$f(x_j)=h(x_j)$ for each $j=1,...,n$. From the definition of
Tychonoff product topology on $G^X$ it follows, that $f\in
cl_{G^X}C_p(X,G)$.

\par {\bf Proposition 2.3.} {\it Let $X$ be a completely $G$-regular
space, $Y\subset X$, $\xi _Y: C_p(X,G)\to C_p(Y,G)$ be a mapping of
the restriction, $\xi _Y(f)=f|_Y$ for each $f\in C_p(X,G)$. Then the
mapping $\xi =\xi _Y$ is continuous and $cl_{C_p(X,G)}\xi
(C_p(X,G))=C_p(Y,G)$.}
\par {\bf Proof.} From definition 2.3. it follows, that
the mapping $\xi $ is continuous. For an arbitrarily given mapping
$g\in C_p(Y,G)$ choose an open neighborhood $W(y_1,...,y_n;U)$ of
the mapping $g$, where $y_1,...,y_n$ belong to $Y$, $U\in {\sf
B}(G)$, $n\in \mathbf{N}$.  Since $X$ is completely $G$-regular
space, then there exists $f\in C_p(X,G)$ such that $f(y_i)=g(y_i)$
for each $i=1,...,n$. For $n=1$ this is evident. For $1<n$ for each
$1\le i\le n$ the subset $P_i := \{ y_j: 1\le j\le n, j\ne i \} $ is
closed in $X$, and for each $q_i$ and $r_i$ in $G$ there exists a
mapping $f_i\in C_p(X,G)$ such that $f_i(y_i)=r_i$, $f_i( P_i)= \{
q_i \} $. Utilizing the ordered products of the form
$d_{i+1}=d_ib_{i+1}$ for $i=1,...,n-1$, where $d_1=b_1$, and
choosing suitable $r_i$ and $q_i$, we deduce the mapping $f$, since
$G$ is the quasigroup. Thus, $\xi (f)\in W(y_1,...,y_n;U)$.

\par {\bf Lemma 2.3.} {\it Let $X$ be a (weakly) completely $G$-regular
space, $A$ be a closed subset in $X$, $B$ be a compact subset in
$X$, $A\cap B=\emptyset $. Then (there exist $q\ne r$ in $G$
correspondingly) for each $q\ne r$ in $G$ there exists a continuous
mapping $f: X\to G$ such that  $f(A)= \{ q \} $ and $f(B)= \{ r \}
$.}
\par {\bf Proof.} Consider an equivalence relation $E$
on $X$ such that $xEy$, if either ($x\in A$ and $y\in A$), or ($x\in
X-A$ and $x=y$). Then there exists the quotient space $Y=X/E$ and
the natural quotient mapping $w: X\to Y$ (see, for example, section
2.4 in \cite{eng}). By virtue of proposition 2.4.3 in \cite{eng} $Y$
is the $T_1$ space, since $A$ is closed in $X$, $w(x)=x$ for each
$x\in X\setminus A$. From proposition 2.4.2 in \cite{eng} it
follows, that for any different points $x\in Y- \{ a \} $ and $y\in
Y$, $x\ne y$, there exists a continuous mapping $f: Y\to G$ such
that $f(x)=r_x$, $f(y)=q_y$, $q_y\ne r_x$ in $G$, where $a=w(A)$.
Thus, the family $C(Y,G)$ of all continuous mappings $f: Y\to G$
separates points in $Y$. From theorem 2.3.20 in \cite{eng} it
follows, that the diagonal $\eta =\Delta _{f\in C(Y,G)}f$ is the
injective continuous mapping $\eta : Y\to G^{C(Y,G)}$. On the other
side, $B$ is compact, consequently, $w(B)$ is compact, hence $w(B)$
is closed in $Y$ according to theorems 3.1.8 and 3.1.10 in
\cite{eng}. Notice that $w(B)\cap w(A)=\emptyset $. Then $\eta
(w(B))\cap \eta (w(A))=\emptyset $, $\eta (w(B))$ is closed in
$G^{C(Y,G)}$, $\eta (w(A))=\eta (a)$. From theorem 2.4 it follows,
that (there exist $q\ne r$ in $G$ correspondingly) for each $q\ne r$
in $G$ there exists the continuous mapping $g: G^{C(Y,G)}\to G$ such
that $g(\eta (w(B)))= \{ r \} $ and $g(\eta (a))=q$, consequently,
$f=g\circ \eta \circ w$ is the continuous mapping $f: X\to G$ such
that $f(A)= \{ q \} $, $f(B)= \{ r \} $.

\par {\bf Proposition 2.4.} {\it Let $X$ be a completely $G$-regular
space, $Y\subset X$, $\xi _Y: C_p(X,G)\to C_p(Y,G)$ be a mapping of
the restriction, and $Y$ be closed in $X$. Then $\xi =\xi _Y$ is the
open mapping from $C_p(X,G)$ onto the subspace $\xi (C_p(X,G))$ in
the space $C_p(Y,G)$. }
\par {\bf Proof.} By virtue of theorem 1.1 it is sufficient to consider the case of
the unital quasigroup $G$. For any different $x_1,...,x_l$ in $Y$
and $x_{l+1},...,x_n$ in $X-Y$ with $n\in \mathbf{N}$ and $f\in
C_p(X,G)$ we consider the open neighborhood $W(x_1,...,x_n;U)$ of
$f$, where $U\in {\sf B}(G)$, $0\le l\le n$. Then $\xi
(W(x_1,...,x_n;U))\subset W(x_1,...,x_l;U)\cap \xi (C_p(X,G))$.
Choose an arbitrary mapping $g\in \xi (C_p(X,G))$ such that $g\in
W(x_1,...,x_l;U)$. Let $h\in C_p(X,G)$ and $\xi (h)=g$. From the
complete $G$-regularity of the space $X$ and the closeness of $Y$ in
$X$ it follows, that there exists a mapping $s\in C_p(X,G)$ such
that $s(Y)= \{ e \} $ and $s(x_j)h(x_j)=f(x_j)$ for each
$j=l+1,...,n$. Then $sh\in W(x_1,...,x_n;U)$ and $\xi (sh)=g$,
consequently, $\xi (W(x_1,...,x_n;U))= W(x_1,...,x_l;U)\cap \xi
(C_p(X,G))$. Thus, the mapping $\xi $ is open.

\par {\bf Proposition 2.5.}
{\it Let the conditions of proposition 2.3 be satisfied, and $Y$ be
everywhere dense in $X$. Then $\xi : C_p(Y,G)\to \xi
(C_p(X,G))=:C_p(Y|X,G)$ is a continuous bijection, where $\xi =\xi
_Y$.}
\par {\bf Proof.} Since $Y$ everywhere dense in $X$, then from
$f\ne g$ in $C_p(X,G)$ it follows, that $f|_Y\ne g|_Y$,
consequently, $\xi (f)\ne \xi (g)$. By virtue of proposition 2.3 the
mapping $\xi =\xi _Y$ is continuous. Thus, the mapping $\xi $ is
continuous, injective and surjective on $\xi (C_p(X,G))$, that is
$\xi $ is the continuous bijection of $C_p(X,G)$ on $C_p(Y|X,G)$.

\par {\bf Corollary 2.6.} {\it If the conditions of proposition 2.5 are satisfied and
$S\subseteq C_p(X,G)$, then the mapping $\xi |_S: S\to \xi
(S)\subseteq C_p(Y|X,G)$ is a homeomorphism.}

\par {\bf Example
2.3.} Let $\tau \ge \aleph _0$, $S$ be a discrete space, $card
(S)=\tau $. From the isomorphism of the topological quasigroups
$(G^{\tau })^{\tau }$ and $G^{\tau }$ it follows, that the
topological quasigroups $C_p(X,G^{\tau })$ and $(C_p(X,G^{\tau
}))^{\tau }$ are isomorphic, since $C_p(X,(G^{\tau })^{\tau })$ is
isomorphic with $C_p(X,G^{\tau })$. Notice that $(C_p(X,G^{\tau
}))^{\tau }$ is isomorphic with $C_p(X\times S, G^{\tau })$. Thus,
from the isomorphism of the topological quasigroups $C_p(X,G^{\tau
})$ and $C_p(Y,G^{\tau })$ generally it does not follow, that the
topological spaces $X$ and $Y$ are homeomorphic.

\par {\bf Definition 2.4.} For a mapping $f$ from a set $X$ into
a set $Y$ let a mapping be given $f^*: G^Y\to G^X$ such that
$f^*(v)(x)=v(f(x))$ for each $v\in G^Y$ and $x\in X$. Then $f^*$ is
called a dual mapping to $f$.

\par {\bf Proposition 2.6.} {\it (a). The dual mapping $f^*$ is continuous.
(b). Let $f(X)=Y$. Then $f^*$ is the homeomorphism from $G^Y$ onto
$f^*(G^Y)$, and $f^*(G^Y)$ is closed in $G^X$.}
\par {\bf Proof.} It evidently follows from definition 2.3
and remark 2.1, since for the discrete topological space $X$ the
spaces $C_p(X,G)$ and $G^X$ are homeomorphic.

\par {\bf Proposition 2.7.} {\it Let $X$, $Y$, $Z$ be completely $G$-regular
spaces, $f: X\to Y$ and $g: X\to Z$ be continuous mappings such that
$f(X)=Y$ and $g(X)=Z$. Then the following conditions are equivalent
(a). $f^*(C(Y,G))\subseteq g^*(C(Z,G))$; (b). There exists a
continuous mapping $h: Z\to Y$ such that $f=h\circ g$.}
\par {\bf Proof.} It is accomplished using definition 2.3
and theorem 2.3 .
\par {\bf Corollary 2.7.} {\it Assume that $X$ and $Y$ are completely $G$-regular spaces,
$f: X\to Y$ is a surjective mapping. Then
\par (a) $f$ is continuous if and only if
$f^*(C(Y,G))\subseteq C(X,G)$; \par (b) if $f$ is the quotient
mapping, then $f^*(C_p(Y,G))$ is the closed subspace in $C_p(X,G)$;
\par (c) $f$ is the continuous bijection if and only if
$f^*(C_p(Y,G))$ is everywhere dense in $C_p(X,G)$;
\par (d) $f$ is a homeomorphism if and only if
$f^*(C_p(Y,G))=C_p(X,G)$. }
\par {\bf Proof.} (a), (b), (d) can be deduced using
proposition 2.7. \par (c). Let $f$ be a continuous bijection, $h\in
C(X,G)$. Choose an open neighborhood $W(x_1,...,x_n;V)$ of the
mapping $h$, where $x_1,...,x_n$ belong to $X$, $n\in \mathbf{N}$,
$V\in {\sf B}(G)$. Using the bijective mapping $f$ we infer, that
$y_i\ne y_j$ for each $i\ne j$, where $f(x_i)=y_i$, $i\in \{ 1,...,
n \} $. Since the space $Y$ is completely $G$-regular, then there
exists a mapping $g$ such that $g(y_i)=h(x_i)$ for each $i\in \{
1,..., n \} $. Therefore $f^*(g)\in W(x_1,...,x_n;V)$, consequently,
$f^*(C_p(Y,G))$ is everywhere dense in $C_p(X,G)$.
\par Vice versa, let $f^*(C_p(Y,G))$ be everywhere dense in $C_p(X,G)$. By
virtue of (a) $f$ is continuous. Assume, that there exist $x_1\ne
x_2$ in $X$ such that $f(x_1)=f(x_2)$. By virtue of theorem 1.1 it
is sufficient to consider the case of the unital quasigroup $G$.
Choose an arbitrary open neighborhood $V$ of the unit element $e$ in
$G$ such that $V\ne G$. There exists an open neighborhood $V_1$ of
the unit element $e$ in $G$ such that $V_1V_1\cup V_1\setminus V_1
\cup V_1/V_1\subset V$, since $G$ is the infinite topological unital
quasigroup being the $T_1$ space (see definition 2.1). Then for each
mapping $h\in f^*(C_p(Y,G))$ we infer that $h(x_1)=h(x_2)$. On the
other side, the space $X$ is completely $G$-regular, consequently,
there exists $g\in C(X,G)$ such that $g(x_1)=e$, $g(x_2)\notin V$.
Then $(W(x_1,x_2;V_1)g)\cap f^*(C_p(Y,G))=\emptyset $. This leads to
the contradiction. Thus, the mapping $f$ is bijective and
continuous, $f(X)=Y$, that is $f$ is the continuous bijection.

\par {\bf Definition 2.5.} Let $X$ be a completely $G$-regular space,
$f: X\to Y$ be a mapping from $X$ onto $Y$. The strongest (i.e.
finest) of all completely $G$-regular topologies $\mathcal{
T}=\mathcal{ T}_{X,f,Y}$ on $Y$, relative to which $f$ is
continuous, is called a $G$-quotient topology on the set $Y$
generated by the mapping $f$. Then the mapping $f$ is called
$G$-quotient, if the topology on $Y$ coincides with the $G$-quotient
topology generated by $f$.
\par {\bf Proposition 2.8.} {\it (a). The family ${\sf B}= \{ h^{-1}(U): h\in
\mathcal{ F}_G, U \mbox{ open in } G \} $ forms a subbase of a
topology $\mathcal{ T}_{X,f,Y}$ (see definition 2.5), where
$\mathcal{ F}_G = \{ h\in G^Y: f^*(h)\in C(X,G) \} $.
\par (b). $\mathcal{ T}_{X,f,Y}$ is the weakest topology on $Y$, for which
all mappings $h$ of $\mathcal{ F}_G$ are continuous. \par (c). If
$\mathcal{ T}^1_Y$ is a topology on $Y$, a mapping $f: X\to
(X,\mathcal{ T}^1_Y)$ is continuous, ${\mathcal F}_G\subset
C((Y,\mathcal{ T}^1_Y),G)$, then $\mathcal{ F}_G = C((Y,\mathcal{
T}^1_Y),G)$.
\par (d). $\mathcal{ T}_{X,f,Y}$ is a unique completely $G$-regular
topology on $Y$, for which $C((Y,\mathcal{ T}_{X,f,Y}),G)=\mathcal{
F}_G$.}
\par {\bf Proof.} This follows from theorems 2.2 and 2.4, and
definition 2.5.
\par {\bf Definition 2.6.} A mapping $f: X\to Y$
of topological completely $G$-regular spaces $X$ and $Y$, where
$f(X)=Y$, is called $G$-functionally closed, if $f^*(C(Y,G))$ is
closed in $C_p(X,G)$.
\par {\bf Proposition 2.9.} {\it Let $X$ and $Y$ be completely $G$-regular spaces.
A mapping $f$ from $X$ onto $Y$ is a $G$-quotient mapping if and
only if $f$ is $G$-functionally closed.}
\par {\bf Proof.} Let a mapping $f$ be $G$-functionally closed.
By virtue of proposition 2.2 $C(Y,G)$ is everywhere dense in $G^Y$.
According to corollary 2.4(a) the mapping $f$ is continuous.
Therefore $f^*(C(Y,G))$ is everywhere dense in $f^*(G^Y)$. On the
other side, $f^*(C(Y,G))\subseteq E_G$, where $E_G=C(X,G)\cap
f^*(G^Y)$, since $f^*(C(Y,G))\subseteq C(X,G)$. Thus, $f^*(C(Y,G))$
is everywhere dense in $E_G$ relative to the topology inherited from
$C_p(X,G)$. Since the mapping $f$ is $G$-functionally closed, then
$f^*(C(Y,G))$ is closed in $C_p(X,G)$, consequently,
$f^*(C(Y,G))=E_G$. From the injectivity of the mapping $f^*$ we
deduce that $C(Y,G)=\mathcal{ F}_G$. The space $Y$ is completely
$G$-regular. Therefore from proposition 2.8(d) it follows, that the
topological space $Y$ has the $G$-quotient topology, generated by
the mapping $f$ (see also definition 2.5). Suppose now that the
mapping $f$ is $G$-quotient. By virtue of proposition 2.8(d)
$C(Y,G)=\mathcal{ F}_G$, consequently, $f^*(C(Y,G))=f^*(\mathcal{
F}_G)$ is the closed subset in $C_p(X,G)$. Thus, $f$ is the
$G$-functionally closed mapping.

\par {\bf Remark 2.3.} Let $G$ be a topological
$T_1$ quasigroup, $X$ be a set, $v_x(f)=f(x)$ for each $x\in X$,
$f\in G^X$, $\mathcal{ F}\subseteq G^X$. Then the mapping $v_x:
\mathcal{ F}\to G$ is continuous for each given $x$ according to
Tychonoff product topology on $G^X$. There also exists a canonical
evaluation mapping $\psi _{\mathcal F}: X\to G^{\mathcal F}$ such
that $\psi _{\mathcal F}(x)=v_x$ for each $x\in X$. If ${\mathcal
F}$ is given, then for brevity it is possible to write $\psi $
instead of $\psi _{\mathcal F}$.

\par {\bf Proposition 2.10.} {\it If $X$ is a
completely $G$-regular space, ${\mathcal F}\subseteq C_p(X,G)$, then
the mapping $\psi _{\mathcal F}: X\to C_p({\mathcal F},G)$ is
continuous.}
\par {\bf Proof.} By virtue of theorem 1.1 it is sufficient to consider
the case of the unital quasigroup $G$. For each $x$ in $X$ and an
open neighborhood $V$ of the unit element $e$ in $G$ such that $V\ne
G$, $f_1,...,f_n$ belong $C_p(X,G)$ there exists an open
neighborhood $W(v_x;f_1,...,f_n;V,G) = \{ h\in C({\mathcal F},G):
\forall j=1,...,n; h(f_j)\in v_x(f_j)V \} $ of $v_x$ in
$C_p({\mathcal F},G)$, where $n\in {\mathbf N}$. Notice that
$\bigcap_{j=1}^nf_j^{-1}(f_j(x)V)=\psi _{\mathcal
F}^{-1}(W(v_x;f_1,...,f_n;V,G))$ is open in $X$ and contains $x$,
since $af\in C_p(X,G)$ and $fs\in C_p(X,G)$ for each $a\in G$, $f$
and $s$ in $C_p(X,G)$. Therefore the mapping $\psi _{\mathcal F}$ is
continuous.

\par {\bf Definition 2.7.} Let $X$ be $Y$ sets, ${\mathcal F}$ be a
family of mappings from $X$ into $Y$. The family ${\mathcal F}$ is
called separating, if for each different points $x_1\ne x_2$ in $X$
there exists $f\in {\mathcal F}$ such that $f(x_1)\ne f(x_2)$. If
$X$ and $Y$ are topological spaces and for each $x\in X$, $A\subset
X$ with $x\notin cl_XA$ there exists $f\in {\mathcal F}$ such that
$f(x)\notin cl_Yf(A)$, then the family ${\mathcal F}$ is called
regular. Notice that, if ${\mathcal F}\subseteq C(X,G)$ and $\psi
_{\mathcal F}: X\to \psi _{\mathcal F}(X)\subseteq C_p({\mathcal
F},G)$ is a homeomorphism, then the family of $G$-valued mappings
${\mathcal F}$ is called generating. If $B$ is a subspace in $X$ and
each continuous mapping $f: B\to G$ has a continuous extension $f_X:
X\to G$, then $B$ is called a $(C,G)$-embedding into $X$.

\par {\bf Remark 2.4.} From definition 2.7 and remark 2.3
it follows, that $\psi _{\mathcal F}(X)$ is the separating family
$G$-valued mappings on ${\mathcal F}$. \par {\bf Proposition 2.11.}
{\it Assume that  $X$ is a completely $G$-regular space, ${\mathcal
F}\subseteq C(X,G)$. Then
\par (a). if ${\mathcal F}$ is a separating family, $\psi _{\mathcal F}(X)$ is an embedding into
$C_p({\mathcal F},G)$, then $\psi _{\mathcal F}: X\to \psi
_{\mathcal F}(X)$ is bijective and continuous.
\par (b). If ${\mathcal F}$ is a regular family, $\psi _{\mathcal
F}(X)$ is an embedding into $C_p({\mathcal F},G)$, then $\psi
_{\mathcal F}: X\to \psi _{\mathcal F}(X)$ is a homeomorphism.}
\par {\bf Proof.} (a). The family ${\mathcal F}$ is separating,
consequently, for each $x_1\ne x_2$ in $X$ there exists $f\in
{\mathcal F}$ with $f(x_1)\ne f(x_2)$; consequently, $v_{x_1}(f)\ne
v_{x_2}(f)$. That is $v_{x_1}\ne v_{x_2}$. \par (b). There exists an
inverse mapping $\psi _{\mathcal F}^{-1}: \psi _{\mathcal F}(X)\to
X$, since $\psi _{\mathcal F}: X\to \psi _{\mathcal F}(X)$ is
bijective according to (a). From the regularity of the family
${\mathcal F}$ it follows, that for each $x\in X$ and $U$ open in
$X$ with $x\in U$, $U\ne X$, there exist an open subset $V$ in $G$
and a mapping $f\in {\mathcal F}$ such that $(f(x)V)\cap
(f(X-U))=\emptyset $. For an open subset $W=W(v_x;f;V,G)=\{ h\in
C({\mathcal F},G): h(f)\in v_x(f)V \} $ in $C_p({\mathcal F},G)$ we
infer that $v_x\in W$ and $\psi _{\mathcal F}^{-1}(W\cap \psi
_{\mathcal F}(X))\subseteq U$, consequently, $\psi _{\mathcal
F}^{-1}$ is the continuous mapping. Applying proposition 2.10 and
the assertion (a) of this proposition we deduce, that $\psi
_{\mathcal F}: X\to \psi _{\mathcal F}(X)\subseteq C_p({\mathcal
F},G)$ is the homeomorphism. \par {\bf Corollary 2.8.} {\it If a
space $X$ is completely $G$-regular, ${\mathcal F}=C_p(X,G)$, then
$X$ is homeomorphic to $\psi (X)$, where $\psi =\psi _{\mathcal F}$,
$\psi (X)$ is the embedding into $C_p(C_p(X,G),G)$.}
\par {\bf Proof.} It follows from proposition 2.11 and remark 2.3.

\par {\bf Remark 2.5.} It is worth to mention that theorems 2.2 - 2.8,
propositions 2.1 - 2.11, lemmas 2.1, 2.2, corollaries 2.1, 2.3 - 2.8
can be transferred on the case of the topological $T_1$ unital left
quasigroup possessing a uniformity compatible with its topology.

\subsection{Conclusion.} All main results of this article are obtained for the first time.
They can be used for further studies of topological spaces
\cite{adbtb,eng}, topological quasigroups \cite{mntpscfb},
noncommutative analysis \cite{nai,reedsim,ludqimqgsmj23,lusmfn2006},
operator theory \cite{ludseabmmiop15}, microbundles
\cite{ludkmtrta19}, representations of topological quasigroups
\cite{nai,hew,sabininb99,lusmfn2006}, in noncommutative geometry
\cite{lawmich,sabininb99}, PDEs \cite{ludhpdecv13,ludspdecv16},
nonassociative quantum field theory
\cite{dnasbdjglta07,dtmnaqmahep07}, nonassociative quantum mechanics
and quantum gravity \cite{mssjmp14,knasqmcstjmp99}, gauge theory and
the great unification theory
\cite{guertzeb,hdtjmp91,mgtnasjmp05,cgostjmp07,incsybea21,nichitaax19,ninybeut14},
in informatics and coding theory
\cite{blautrctb,gkmngcnagdm04,mmnnascdjms20,pualdb94,srwseabm14},
since they are based on algebraic topology of binary systems (see
also introduction).

\end{document}